\documentclass[a4paper, 9pt, conference]{ieeeconf}

\usepackage{amsmath}
\usepackage{accents}
\usepackage{amssymb}
\usepackage{graphicx}
\usepackage{setspace}
\usepackage{mathrsfs}

\newcommand{\ud}{\mathrm{d}}

\newcommand{\bfa}{\mathbf{a}}
\newcommand{\bfb}{\mathbf{b}}

\newcommand{\bff}{\mathbf{f}}

\newcommand{\bfx}{\mathbf{x}}

\newcommand{\bfA}{\mathbf{A}}

\newcommand{\bfF}{\mathbf{F}}
\newcommand{\bfG}{\mathbf{G}}

\newcommand{\bfI}{\mathbf{I}}

\newcommand{\bfM}{\mathbf{M}}
\newcommand{\bfP}{\mathbf{P}}
\newcommand{\bfQ}{\mathbf{Q}}

\newcommand{\bfTh}{\mathbf{\Theta}}
\newcommand{\zeros}{\mathbf{0}}

\newcommand{\dm}[1]{d_{\mathbf{M}(#1)}}
\newcommand{\tr}{\mathrm{tr}}

\newcommand{\expect}{\mathbb{E}}
\newcommand{\expectx}{\mathbb{E}_{\bfx}}

\makeatletter
\renewcommand\theequation{\thesection.\arabic{equation}}
\@addtoreset{equation}{section}
\makeatother
\newtheorem{prop}{Proposition}
\newtheorem{theorem}{Theorem}

\title{Stochastic Contraction in Riemannian Metrics}

\author{
Quang-Cuong Pham\\
{\small University of Tokyo}\\
\and Jean-Jacques Slotine\\
{\small Massachusetts Institute of Technology}\\
}

\begin{document}

\maketitle

\begin{abstract}
  Stochastic contraction analysis is a recently developed tool for
  studying the global stability properties of nonlinear stochastic
  systems, based on a differential analysis of convergence in an
  appropriate metric.  To date, stochastic contraction results and
  sharp associated performance bounds have been established only in
  the specialized context of state-independent metrics, which
  restricts their applicability.  This paper extends stochastic
  contraction analysis to the case of general time- and
  state-dependent Riemannian metrics, in both discrete-time and
  continuous-time settings, thus extending its applicability to a
  significantly wider range of nonlinear stochastic dynamics.
\end{abstract}

\section{Introduction}

Contraction theory provides a body of analytical tools to study the
stability and convergence of nonlinear dynamical
systems~\cite{LS98automatica}. Based on a differential analysis of
convergence, it allows global stability properties of a nonlinear
system to be concluded from the system's linearization at all points
in some appropriate metric. Historically, basic convergence results on
contracting systems can be traced back to the numerical analysis
literature~\cite{Lew49ajm,Har64book,Dem61smm}.  Recently, contraction
theory has been extended to \emph{stochastic} differential
systems~\cite{PhaX09tac}. This development has led to a number of
practically important applications, such as the design of observers
for nonlinear stochastic systems~\cite{DanX13tac}, or the study of
synchronization in networks of noisy oscillators~\cite{TabX10pcb}.

The stochastic contraction theorems have been formulated so far in the
specialized context of \emph{state-independent}
metrics~\cite{PhaX09tac}. Yet, more general \emph{state-dependent}
Riemannian metrics can be central to some systems, and in fact the
original deterministic contraction theorems were derived in this
general context~\cite{LS98automatica}. Some practical nonlinear
dynamics can be most easily studied by choosing appropriate
state-dependent metrics (cf. e.g.~\cite{AR03tac,DanX13tac}), and from
a theoretical perspective, the contraction properties of some systems
can \emph{only} be observed in a state-dependent Riemannian
metric~\cite{SP12tac}.

Recently, an attempt has been made to extend the stochastic
contraction results of~\cite{PhaX09tac} to state-dependent
metrics~\cite{DanX13tac}. However, since in the estimation of the
distance between two trajectories the derivation did not consider
geodesics between these trajectories but instead used straight lines,
the bounds obtained are not ``optimal'' (in a sense made precise in
Remark~3.3). Here, we prove the stochastic contraction theorems in the
case of general time- and state-dependent Riemannian metrics by
studying the evolution of the geodesics under the combined effects of
the noise and the contracting flow, which allows ``optimal'' bounds to
be obtained.

In section~\ref{sec:discrete}, we study the contraction properties of
discrete-time stochastic difference systems. Then, in
section~\ref{sec:continuous}, we address the case of continuous-time It\^o
stochastic differential systems by using a discrete/continuous
limiting argument. Finally, section~\ref{sec:conclusion} offers
brief concluding remarks.

\section{Discrete stochastic contraction}
\label{sec:discrete}

We first state and prove a proposition (see also~\cite{AR03tac}), which
makes explicit the original deterministic discrete contraction theorem
(see section 5 of~\cite{LS98automatica}).

%%%%%%%%%%%%%%%%%%%%%%%%%%%%%%%%%%%%%%%%%%%%%%%%%%%%%%%%%%
%%%%%%%%%%%%%%%%%%%%%% Prop 1 %%%%%%%%%%%%%%%%%%%%%%%%%%%%
%%%%%%%%%%%%%%%%%%%%%%%%%%%%%%%%%%%%%%%%%%%%%%%%%%%%%%%%%%

\begin{prop}[and definition]
  \label{prop:ineq}
  Consider two uniformly positive definite metrics
  $\bfM_i=\bfTh_i^\top\bfTh_i$ ($i=1,2$) defined over $\mathbb{R}^n$ and
  a smooth function $\bff:\mathbb{R}^n \to\mathbb{R}^n$. The
  \emph{generalized Jacobian} of $\bff$ in the metrics
  $(\bfM_1,\bfM_2)$ is defined by
  \[
  \bfF=\bfTh_2\frac{\partial \bff}{\partial \bfa}\bfTh_1^{-1}.
  \]

  Assume now that $\bff$ is \emph{contracting} in the metrics
  $(\bfM_1,\bfM_2)$ with rate $\mu$, i.e.
  \[
  \forall \bfa\in\mathbb{R}^n 
  \quad \lambda_{\max}(\bfF(\bfa)^\top\bfF(\bfa)) \leq \mu,
  \]  
  where $\lambda_{\max}(\bfA)$ denotes the largest eigenvalue of a
  given matrix $\bfA$.  Then for all $\bfa,\bfb\in\mathbb{R}^n$, one
  has
  \[
  d^2_{\bfM_2}(\bff(\bfa),\bff(\bfb)) \leq \mu d^2_{\bfM_1}(\bfa,\bfb),
  \]
  where $d_{\bfM}$ denotes the distance associated with the metric
  $\bfM$.
\end{prop}

%%%%%%%%%%%%%%%%%%%%%%%%%%% Proof %%%%%%%%%%%%%%%%%%%%%%%%%%%%%%%%%

\paragraph*{Proof} Since $\bfM_1$ is uniformly positive definite, there
exists a $C^1$-continuous curve (a geodesic) $\Gamma : [0,1]\to
\mathbb{R}^n$ such that $\Gamma(0)=\bfa$ and $\Gamma(1)=\bfb$ and
\[
d^2_{\bfM_1}(\bfa,\bfb)=\int_0^1 \left(\frac{\partial \Gamma}{\partial
    u}(u)\right)^\top \bfM_1(\Gamma(u))\left(\frac{\partial \Gamma}{\partial
    u}(u)\right) du.
\]
Next, since $\bff$ is a smooth function, $\bff(\Gamma)$ is also a
$C^1$-continuous curve. By the definition of the distance, one then
has
\[
d^2_{\bfM_2}(\bff(\bfa),\bff(\bfb))\leq \int_0^1 \left(\frac{\partial
      \bff(\Gamma)}{\partial u}(u)\right)^\top \bfM_2(\bff(\Gamma(u)))\left(\frac{\partial
      \bff(\Gamma)}{\partial u}(u)\right) du.
\]

Remark on the other hand that,  by the chain rule,
\[
\frac{\partial \bff(\Gamma)}{\partial u}(u)=\frac{\partial
  \bff}{\partial \bfa} \frac{\partial \Gamma}{\partial u}(u),
\]
which leads to
% XXXX
% \[
%   \begin{array}{rcl}
% \int_0^1 \left(\frac{\partial
%       \bff(\Gamma)}{\partial u}(u)\right)^\top \bfM_2\left(\frac{\partial
%       \bff(\Gamma)}{\partial u}(u)\right) du
%     &=&\int_0^1 \left(
%         \frac{\partial \Gamma}{\partial u}^\top
%         \frac{\partial \bff}{\partial \bfa}^\top\bfTh_2^\top\bfTh_2
%         \frac{\partial \bff}{\partial \bfa} 
%         \frac{\partial \Gamma}{\partial u} \right) du \\
%     &=&\int_0^1
%       \left(\frac{\partial \Gamma}{\partial u}^\top\bfTh_1^\top\right)
%       \bfF^\top\bfF
%       \left(\bfTh_1\frac{\partial \Gamma}{\partial u}\right) du \\
%     &\leq& \int_0^1 \mu \left (\frac{\partial \Gamma}{\partial u}^\top 
%         \bfTh_1^\top\bfTh_1
%         \frac{\partial \Gamma}{\partial u} \right) du \\
%     &=&\mu d^2_{\bfM_1}(\bfa,\bfb) \quad \Box
%   \end{array}
% \]
\[
  \begin{array}{rcl}
&&\int_0^1 \left(\frac{\partial
      \bff(\Gamma)}{\partial u}(u)\right)^\top \bfM_2\left(\frac{\partial
      \bff(\Gamma)}{\partial u}(u)\right) du\\
    &&=\int_0^1 \left(
        \frac{\partial \Gamma}{\partial u}^\top
        \frac{\partial \bff}{\partial \bfa}^\top\bfTh_2^\top\bfTh_2
        \frac{\partial \bff}{\partial \bfa} 
        \frac{\partial \Gamma}{\partial u} \right) du \\
    &&=\int_0^1
      \left(\frac{\partial \Gamma}{\partial u}^\top\bfTh_1^\top\right)
      \bfF^\top\bfF
      \left(\bfTh_1\frac{\partial \Gamma}{\partial u}\right) du \\
    &&\leq \int_0^1 \mu \left (\frac{\partial \Gamma}{\partial u}^\top 
        \bfTh_1^\top\bfTh_1
        \frac{\partial \Gamma}{\partial u} \right) du \\
    &&=\mu d^2_{\bfM_1}(\bfa,\bfb) \quad \Box
  \end{array}
\]

We now state and prove a proposition which relates metrics and noise.

%%%%%%%%%%%%%%%%%%%%%%%%%%%%%%%%%%%%%%%%%%%%%%%%%%%%%%%%%%
%%%%%%%%%%%%%%%%%%%%%% Prop 2 %%%%%%%%%%%%%%%%%%%%%%%%%%%%
%%%%%%%%%%%%%%%%%%%%%%%%%%%%%%%%%%%%%%%%%%%%%%%%%%%%%%%%%%

\begin{prop}
  \label{prop:eta}
  Consider a uniformly positive definite metric $\bfM$ defined over
  $\mathbb{R}^n$. Let $\sigma$ be a matrix-valued function
  $\mathbb{R}^n\rightarrow\mathbb{R}^{nd}$, $\eta_1,\eta_2$ two
  independent $d$-dimensional Gaussian random variables with
  $\eta_i\sim \mathscr{N}(\bf0,\bfI)$, and
  $\bfa,\bfb\in\mathbb{R}^n$. Assume that
  \[
  \forall \bfa\in\mathbb{R}^n\quad \tr(\sigma(\bfa)^\top\bfM(\bfa)\sigma(\bfa))\leq D,
  \]
  \[
  \textrm{then one has}\quad 
  \expect\left[d^2_\bfM(\bfa+\sigma(\bfa)\eta_1,\bfb+\sigma(\bfb)\eta_2)\right] \leq 
  d^2_\bfM(\bfa,\bfb) + 2D.
  \]
\end{prop}

%%%%%%%%%%%%%%%%%%%%%%%%%%% Proof %%%%%%%%%%%%%%%%%%%%%%%%%%%%%

\paragraph*{Proof} As previously, since $\bfM$ is uniformly positive
definite, there exists a $C^1$-continuous curve $\Gamma : [0,1]\to
\mathbb{R}^n$ such that $\Gamma(0)=\bfa$ and $\Gamma(1)=\bfb$ and
\[
d^2_{\bfM}(\bfa,\bfb)=\int_0^1 \left(\frac{\partial \Gamma}{\partial
    u}\right)^\top \bfM\left(\frac{\partial \Gamma}{\partial
    u}\right) du.
\]

Consider the curve $\Gamma_\eta: [0,1]\to
\mathbb{R}^n$ defined by
\[
\forall u\in[0,1]\quad
\Gamma_\eta(u)=\Gamma(u)+(1-u)\sigma(\bfa)\eta_1+u\sigma(\bfb)\eta_2. 
\]
It is clear that $\Gamma_\eta$ is $C^1$-continuous and verifies
$\Gamma_\eta(0)=\bfa+\sigma(\bfa)\eta_1$ and
$\Gamma_\eta(1)=\bfb+\sigma(\bfb)\eta_2$. Thus, by the definition of
the distance, one has
%XXXX
% \begin{eqnarray}
%   d^2_\bfM(\bfa+\eta_1,\bfb+\eta_2)&\leq& \int_0^1 \left(\frac{\partial \Gamma_\eta}{\partial
%       u}\right)^\top \bfM\left(\frac{\partial \Gamma_\eta}{\partial
%       u}\right) du\nonumber\\
%   &=&\int_0^1 \left(\frac{\partial \Gamma}{\partial
%       u}+(\sigma(\bfb)\eta_2-\sigma(\bfa)\eta_1)\right)^\top
%   \bfM\left(\frac{\partial \Gamma}{\partial 
%       u}+(\sigma(\bfb)\eta_2-\sigma(\bfa)\eta_1)\right) du\nonumber\\
%   &=&d^2_{\bfM}(\bfa,\bfb)
%   +2(\sigma(\bfb)\eta_2-\sigma(\bfa)\eta_1)^\top\int_0^1 \bfM
%   \left(\frac{\partial \Gamma}{\partial 
%       u}\right)du\nonumber\\
%   &&-2(\sigma(\bfb)\eta_2)^\top\left(\int_0^1\bfM
%     du\right)(\sigma(\bfa)\eta_1) \nonumber\\
%   &&+\int_0^1(\sigma(\bfa)\eta_1)^\top\bfM(\sigma(\bfa)\eta_1)
%     du+\int_0^1(\sigma(\bfb)\eta_2)^\top\bfM(\sigma(\bfb)\eta_2)
%     du.\nonumber
% \end{eqnarray}
\[  
d^2_\bfM(\bfa+\eta_1,\bfb+\eta_2)
\]
\begin{eqnarray}
  &\leq& \int_0^1 \left(\frac{\partial \Gamma_\eta}{\partial
      u}\right)^\top \bfM\left(\frac{\partial \Gamma_\eta}{\partial
      u}\right) du\nonumber\\
  &=&\int_0^1 \left(\frac{\partial \Gamma}{\partial
      u}+(\sigma(\bfb)\eta_2-\sigma(\bfa)\eta_1)\right)^\top
  \bfM\nonumber\\
  &&\left(\frac{\partial \Gamma}{\partial 
      u}+(\sigma(\bfb)\eta_2-\sigma(\bfa)\eta_1)\right) du\nonumber\\
  &=&d^2_{\bfM}(\bfa,\bfb)
  +2(\sigma(\bfb)\eta_2-\sigma(\bfa)\eta_1)^\top\int_0^1 \bfM
  \left(\frac{\partial \Gamma}{\partial 
      u}\right)du\nonumber\\
  &&-2(\sigma(\bfb)\eta_2)^\top\left(\int_0^1\bfM
    du\right)(\sigma(\bfa)\eta_1) \nonumber\\
  &&+\int_0^1(\sigma(\bfa)\eta_1)^\top\bfM(\sigma(\bfa)\eta_1)
  du\nonumber\\
  &&+\int_0^1(\sigma(\bfb)\eta_2)^\top\bfM(\sigma(\bfb)\eta_2)
du.\nonumber
\end{eqnarray}

Remark that the second and third terms of the right-hand side vanish
when taking the expectation. As for the fourth and fifth terms, remark
that
\begin{eqnarray}
\label{eq:trace}
(\sigma(\bfa)\eta_1)^\top\bfM(\sigma(\bfa)\eta_1)&=&
\tr\left((\sigma(\bfa)\eta_1)^\top\bfM(\sigma(\bfa)\eta_1)\right)\nonumber\\
&=&
\tr\left(\eta_1^\top\sigma(\bfa)^\top\bfM\sigma(\bfa)\eta_1\right)=\tr(\eta_1^\top\bfQ\eta_1),\nonumber
\end{eqnarray}
where $\bfQ$ is obtained from $\sigma(\bfa)^\top\bfM\sigma(\bfa)$ by an
orthogonal diagonalization. One thus has
\[
\expect\left[\tr(\eta_1^\top\bfQ\eta_1)\right]=\tr(\bfQ)=\tr(\sigma(\bfa)^\top\bfM\sigma(\bfa))\leq D,
\]
which allows to conclude $\Box$

We can now state and prove the discrete stochastic contraction
theorem.

%%%%%%%%%%%%%%%%%%%%%%%%%%%%%%%%%%%%%%%%%%%%%%%%%%%%%%%%%%%%%
%%%%%%%%%%%%%%%%%%%%%% Theorem 1 %%%%%%%%%%%%%%%%%%%%%%%%%%%%
%%%%%%%%%%%%%%%%%%%%%%%%%%%%%%%%%%%%%%%%%%%%%%%%%%%%%%%%%%%%%

\begin{theorem}
  \label{theo:discrete}

  Consider the stochastic difference equation

  \begin{equation}
    \label{eq:main}
    \left\{
    \begin{array}{l}
      \bfa_{k+1}=\bff(\bfa_k,k)+\sigma(\bfa_k,k)w_{k+1}\\
      \bfa_0=\xi, 
    \end{array}
    \right.
  \end{equation}
  where $\bff$ is a $\mathbb{R}^n \times \mathbb{N}\to\mathbb{R}^n$
  function, $\sigma$ is a $\mathbb{R}^n \times
  \mathbb{N}\to\mathbb{R}^{nd}$ matrix-valued function,
  $(w_k)_{k\in\mathbb{N}}$ is a sequence of independent $d$-dimensional
  Gaussian noise vectors, with $w_k\sim \mathscr{N}(\zeros,\bfI)$
  and $\xi$ is a $n$-dimensional random variable independent of the
  $w_k$.
  
  Assume that the system verifies the following two hypotheses:
  
  \begin{description}
  \item[\textbf{(Hd1)}] for all $k\geq 0$, the dynamics $\bff(\bfa,k)$
    is contracting in the metrics $(\bfM_k,\bfM_{k+1})$, with
    contraction rate $\mu$ $(0<\mu<1)$, and the metrics
    $\bfM_k(\bfa)$ are uniformly positive definite in $\bfa$ and
    $k$, with lower bound~$\beta$, i.e.
    \[
    \forall  k\geq 0,\  \bfa\in\mathbb{R}^n\quad 
    \bfa^\top\bfM_k(\bfa)\bfa \geq \beta \|\bfa\|^2;
    \]

  \item[\textbf{(Hd2)}]
    $\tr\left(\sigma(\bfa,k)^\top\bfM(\bfa,k)\sigma(\bfa,k)\right)$ is
    uniformly upper-bounded by a constant~$D$.
 
  \end{description}
  
  Let $(\bfa_k)_{k\in\mathbb{N}}$ and $(\bfb_k)_{k\in\mathbb{N}}$ be
  two trajectories whose initial conditions are given by a probability
  distribution $p(\xi,\xi')$. Then for all $k\geq 0$,
  % 
  % XXXX
  % \begin{eqnarray}
  %   \label{eq:etot}
  %   \expect \left[d^2_{\bfM_k}(\bfa_k,\bfb_k)\right] \leq
  %   \frac{2D}{1-\mu}  + 
  %   \mu^k \int
  %   \left[d^2_{\bfM_0}(\bfa_0,\bfb_0)-\frac{2D}{1-\mu}\right]^+
  %   dp({\bfa_0},{\bfb_0}),
  % \end{eqnarray}
  \begin{eqnarray}
    \label{eq:etot}
    &&\expect \left[d^2_{\bfM_k}(\bfa_k,\bfb_k)\right] \leq
    \frac{2D}{1-\mu}  \nonumber\\
    &+& 
    \mu^k \int
    \left[d^2_{\bfM_0}(\bfa_0,\bfb_0)-\frac{2D}{1-\mu}\right]^+
    dp({\bfa_0},{\bfb_0}),
  \end{eqnarray}
  where $[\cdot]^+=\max(0,\cdot)$. 

  In particular, for all $k \geq 0$,
  \begin{equation}
    \label{eq:e} 
    \expect \left[\|\bfa_k-\bfb_k\|^2 \right] \leq
    \frac{2D}{\beta(1-\mu)} + 
    \frac{\mu^k}{\beta}\expect\left[d^2_{\bfM_0}(\xi,\xi')\right].
  \end{equation}  
\end{theorem}

%%%%%%%%%%%%%%%%%%%%%%%% Proof %%%%%%%%%%%%%%%%%%%%%%%%%%%%%%

\paragraph*{Proof} Taking the conditional expectation given
$(\bfa_0,\bfb_0)=\bfx$ and applying \textbf{(H2d)} and
Proposition~\ref{prop:eta}, one has
%XXXX
% \[
% \begin{array}{rcl}
% \expectx\left[d^2_{\bfM_{k+1}}(\bfa_{k+1},\bfb_{k+1})\right]
% &=&\expectx\left[d^2_{\bfM_{k+1}}\left(\bff(\bfa,k)+\sigma(\bfa,k)w_k,
%           \bff(\bfb,k)+\sigma(\bfb,k)w'_k \right)\right]\\
% &\leq& \expectx\left[d^2_{\bfM_{k+1}}(\bff(\bfa_k),\bff(\bfb_k))\right]+2D,
% \end{array}
% \]
\[
\begin{array}{rcl}
\expectx\left[d^2_{\bfM_{k+1}}(\bfa_{k+1},\bfb_{k+1})\right]
&=&\expectx\left[d^2_{\bfM_{k+1}}\left(\bff(\bfa,k)+\sigma(\bfa,k)w_k,\right.\right.\\
&&\left.\left.          \bff(\bfb,k)+\sigma(\bfb,k)w'_k \right)\right]\\
&\leq& \expectx\left[d^2_{\bfM_{k+1}}(\bff(\bfa_k),\bff(\bfb_k))\right]+2D,
\end{array}
\]
where $w'_k$  has the same distribution as $w_k$ but is independent of
the latter.

On the other hand, from \textbf{(Hd1)} and Proposition
\ref{prop:ineq}, one has
\[
\expectx \left[d^2_{\bfM_{k+1}}(\bff(\bfa_k),\bff(\bfb_k))\right] \leq
\mu \expectx \left[d^2_{\bfM_k}(\bfa_k,\bfb_k)\right].
\]

If one now sets $u_k=\expectx\left[d_{\bfM_k}(\bfa_k,\bfb_k)\right]$
then it follows from the above that
\begin{equation}
  \label{eq:u}
  u_{k+1}\leq \mu u_k+2D.
\end{equation}

Define next $v_k=u_k-2D/(1-\mu)$. Then replacing $u_k$ by
$v_k+2D/(1-\mu)$ in (\ref{eq:u}) leads to $v_{k+1}\leq\mu v_k$. This
implies that $\forall k\geq 0, \ v_k\leq v_0\mu^k \leq
[v_0]^+\mu^k$. Replacing $v_k$ by its expression in terms of $u_k$
then yields
\[
\forall k\geq 0\quad u_k\leq
\frac{2D}{1-\mu}+\mu
^k\left[u_0-\frac{2D}{1-\mu}\right]^+.
\]
Integrating the last inequality with respect to $\bfx$ leads to
(\ref{eq:etot}).  Finally, (\ref{eq:e}) follows from (\ref{eq:etot})
by remarking that
\begin{eqnarray}
  \label{eq:trick}
  \int
  \left[d^2_{\bfM_0}(\bfa_0,\bfb_0)-\frac{2D}{1-\mu}\right]^+
  dp(\bfa_0,\bfb_0) \leq \nonumber \\
  \int
  d^2_{\bfM_0}(\bfa_0,\bfb_0) dp(\bfa_0,\bfb_0)
  = \expect\left[d^2_{\bfM_0}(\xi,\xi')\right],
\end{eqnarray}
\begin{equation}
  \label{eq:trick2}
  \textrm{and that}, \ 
\|\bfa_k-\bfb_k\|^2 \leq
\frac{1}{\beta}d_{\bfM_k}^2(\bfa_k,\bfb_k) \quad \Box  
\end{equation}

\paragraph*{Remark 2.1 [Relaxing the uniform bound on the noise]} Assume
that the initial conditions are contained in a region $U$, then
\textbf{(Hd2)} can in fact be replaced by~\cite{Her13perso}
\[
\forall k\geq 0,\ \bfa \in U   \quad \expect\left[
  \tr\left(\sigma(\bfa_k,k)^\top\bfM_k(\bfa_k)\sigma(\bfa_k,k)\right) \ | \
  \bfa_0=\bfa \right] \leq D.
\]

\section{Continuous stochastic contraction}
\label{sec:continuous}

Based on the discrete stochastic contraction theorem just established,
we can now state and prove the continuous stochastic contraction
theorem in general Riemannian metrics.

Consider the It\^o stochastic differential equation
\begin{equation}
  \label{eq:main-cont}
  \left\{\begin{array}{l}
      d\bfa=\bff(\bfa,t)dt+\sigma(\bfa,t)dW\\
      \bfa(0)=\xi.
    \end{array}  
  \right.
\end{equation}

To ensure existence and uniqueness of solutions to equation
(\ref{eq:main}), we assume the following standard conditions on $\bff$
and $\sigma$:

\emph{Lipschitz condition:} There exists a constant $K_1>0$ such that
\[ \forall t\geq 0, \ \bfa,\bfb \in\mathbb{R}^n \quad
\|\bff(\bfa,t)-\bff(\bfb,t)\|+\|\sigma(\bfa,t)-\sigma(\bfb,t)\| \leq
K_1\|\bfa-\bfb\|;
\]

\emph{Restriction on growth:} There exists a constant $K_2>0$ such that
\[ \forall t\geq 0, \ \bfa \in\mathbb{R}^n \quad
\|\bff(\bfa,t)\|^2+\|\sigma(\bfa,t)\|^2 \leq K_2(1+\|\bfa\|^2).
\]

%%%%%%%%%%%%%%%%%%%%%%%%%%%%%%%%%%%%%%%%%%%%%%%%%%%%%%%%%%%%%
%%%%%%%%%%%%%%%%%%%%%% Theorem 2 %%%%%%%%%%%%%%%%%%%%%%%%%%%%
%%%%%%%%%%%%%%%%%%%%%%%%%%%%%%%%%%%%%%%%%%%%%%%%%%%%%%%%%%%%%

\begin{theorem}
  \label{theo:main-cont}
  Assume that system (\ref{eq:main-cont}) verifies the following two
  hypotheses:

  \begin{description}
  \item[\textbf{(Hc1)}] for all $t\geq 0$, the dynamics $\bff(\bfa,t)$
    is contracting in the time- and state-dependent metric
    $\bfM(\bfa,t)=\bfTh^\top(\bfa,t)\bfTh(\bfa,t)$, with contraction
    rate $\lambda$ $(\lambda>0)$, i.e.
    %XXXX
    % \[
    % \forall t\geq 0 ,\ \bfa\in\mathbb{R}^n \quad \lambda_{\max}\left(\left[
    %   \left(\dot{\bfTh}(\bfa,t)+ \bfTh(\bfa,t)\frac{\partial \bff}{\partial
    %       \bfa} \right)\bfTh(\bfa,t)^{-1} \right]_s\right) \leq -\lambda,
    % \]
    \[
    \forall t\geq 0 ,\ \bfa\in\mathbb{R}^n 
    \]
    \[ 
    \lambda_{\max}\left(\left[
      \left(\dot{\bfTh}(\bfa,t)+ \bfTh(\bfa,t)\frac{\partial \bff}{\partial
          \bfa} \right)\bfTh(\bfa,t)^{-1} \right]_s\right) \leq -\lambda,
    \]
    where $\bfA_s=\frac{1}{2}(\bfA^\top+\bfA)$ denotes the symmetric
    part of a given matrix~$\bfA$.  Furthermore, the metric
    $\bfM(\bfa,t)$ is positive definite uniformly in $\bfa$ and~$t$,
    with lower bound $\beta$;

  \item[\textbf{(Hc2)}]
    $\tr\left(\sigma(\bfa,t)^\top\bfM(\bfa,t)\sigma(\bfa,t)\right)$
    is uniformly upper-bounded by a constant~C.
  \end{description}

  Let $\bfa(t)$ and $\bfb(t)$ be two trajectories whose initial
  conditions are independent of $W$ and given by a probability
  distribution $p(\xi,\xi')$. Then for all $T\geq 0$,
  %
  % XXXX
  % \begin{eqnarray}
  %   \label{eq:gen1}
  %   \expect \left[ \dm{T}^2(\bfa(T),\bfb(T)) \right] \leq 
  %   \frac{C}{\lambda} +
  %    e^{-2\lambda T}\int
  %     \left[\dm{0}^2(\bfa_0,\bfb_0)
  %       -\frac{C}{\lambda}\right]^+dp(\bfa_0,\bfb_0).
  % \end{eqnarray}
  \begin{eqnarray}
    \label{eq:gen1}
    &&\expect \left[ \dm{T}^2(\bfa(T),\bfb(T)) \right] \leq 
    \frac{C}{\lambda} \nonumber\\
    &&+
     e^{-2\lambda T}\int
      \left[\dm{0}^2(\bfa_0,\bfb_0)
        -\frac{C}{\lambda}\right]^+dp(\bfa_0,\bfb_0).
  \end{eqnarray}
   
  In particular, for all $T \geq 0$,
  \begin{equation} 
    \label{eq:gen2}
    \expect \left[ \|\bfa(T)-\bfb(T)\|^2 \right] \leq
    \frac{C}{\beta\lambda} +
    \frac{e^{-2\lambda T}}{\beta}\expect\left[ \dm{0}^2(\xi,\xi') \right].
  \end{equation}
\end{theorem}

%%%%%%%%%%%%%%%%%%%%%%%% Proof %%%%%%%%%%%%%%%%%%%%%%%%%%%%%%

\paragraph*{Proof} Fix $(\bfa(0),\bfb(0))=\bfx\in\mathbb{R}^{2d}$ and
$T\geq 0$. We first discretize the time interval $[0,T]$ into $N$
equal intervals of length $\delta=T/N$ and consider the two sequences
$(\bfa^\delta_k)_{k\in\mathbb{N}}$, $(\bfb^\delta_k)_{k\in\mathbb{N}}$
defined by
%
% XXXX
% \begin{equation}
%   \label{eq:SDEdiscr}
%   \left\{
%   \begin{array}{l}
%     \bfa^\delta_{k+1}=\bfa^\delta_k+\delta \bff(\bfa^\delta_k,k\delta) +
%     \sigma(\bfa^\delta_k,k\delta)w^\delta_k\\
%     \bfa^\delta_0=\bfa(0)
%   \end{array}
%   \right.
%   \left\{
%   \begin{array}{l}
%     \bfb^\delta_{k+1}=\bfb^\delta_k+\delta \bff(\bfb^\delta_k,k\delta) +
%     \sigma(\bfb^\delta_k,k\delta)w'^\delta_k\\
%     \bfb^\delta_0=\bfb(0),
%   \end{array}
%   \right.
% \end{equation}
\begin{eqnarray}
  \label{eq:SDEdiscr}
  \left\{
  \begin{array}{l}
    \bfa^\delta_{k+1}=\bfa^\delta_k+\delta \bff(\bfa^\delta_k,k\delta) +
    \sigma(\bfa^\delta_k,k\delta)w^\delta_k\\
    \bfa^\delta_0=\bfa(0)
  \end{array}
  \right.\nonumber\\
  \left\{
  \begin{array}{l}
    \bfb^\delta_{k+1}=\bfb^\delta_k+\delta \bff(\bfb^\delta_k,k\delta) +
    \sigma(\bfb^\delta_k,k\delta)w'^\delta_k\\
    \bfb^\delta_0=\bfb(0),
  \end{array}
  \right.
\end{eqnarray}
where $(w^\delta_k)_{k\in\mathbb{N}}$ and
$(w'^\delta_k)_{k\in\mathbb{N}}$ are two sequences of random variables
defined by $w^\delta_k=W((k+1)\delta)-W(k\delta)$ and
$w'^\delta_k=W'((k+1)\delta)-W'(k\delta)$. Note that, since $W$ and
$W'$ are two independent Wiener processes, $(w_k)_{k\in\mathbb{N}}$
and $(w'_k)_{k\in\mathbb{N}}$ are two sequences of independent
Gaussian random variables with distribution
$\mathcal{N}(\zeros,\delta\bfI)$. Note also that, by the strong
convergence of the Euler-Maruyama scheme (cf.~\cite{KP92book},
p.~342), one has
\begin{equation}
  \label{eq:lim}
  \lim_{\delta\to 0}\expect_\bfx\left[\|\bfa^\delta_N-\bfa(T)\|^2\right] =0.
\end{equation}

Hypothesis \textbf{(Hc2)} implies that system~(\ref{eq:SDEdiscr})
satisfies \textbf{(Hd2)} with $D=\delta C$.  To verify \textbf{(Hc1)},
denote by $\bfG_k(\bfa)$ the generalized Jacobian matrix
of~(\ref{eq:SDEdiscr}) at step $k$. Denoting $t=k\delta$, one has
%XXXX
% \[
% \bfG_k(\bfa)=\bfTh(\bfa,t+\delta)\frac{\partial (\bfa+\delta \bff(\bfa,t))}{\partial \bfa} 
% \bfTh(\bfa,t)^{-1}=\bfTh(t+\delta)
% \left(\bfI+\delta\frac{\partial \bff}{\partial \bfa}\right)\bfTh(t)^{-1}.
% \]
\[
\bfG_k(\bfa)=\bfTh(\bfa,t+\delta)\frac{\partial (\bfa+\delta \bff(\bfa,t))}{\partial \bfa} 
\bfTh(\bfa,t)^{-1}\]
\[=\bfTh(t+\delta)
\left(\bfI+\delta\frac{\partial \bff}{\partial \bfa}\right)\bfTh(t)^{-1}.
\]
Remark that we have dropped the argument $\bfa$ for convenience. One
can next rewrite $\bfG_k^\top\bfG_k=\bfA_0+\delta \bfA_1$, with
\[
% XXXX
% \begin{array}{rcl}
%   \bfA_0&=&\left(\bfTh(t)^{-1}\right)^\top\bfTh(t+\delta)^\top
%   \bfTh(t+\delta)\bfTh(t)^{-1};\\
%   \bfA_1&=&\delta\left(\bfTh(t)^{-1}\right)^\top\left(
%     \bfTh(t+\delta)^\top\bfTh(t+\delta)\frac{\partial \bff}{\partial \bfa}  +
%     \left(\frac{\partial \bff}{\partial \bfa}\right)^\top\bfTh(t+\delta)^\top\bfTh(t+\delta)
%   \right)\bfTh(t)^{-1}.
% \end{array}   
\begin{array}{rcl}
  \bfA_0&=&\left(\bfTh(t)^{-1}\right)^\top\bfTh(t+\delta)^\top
  \bfTh(t+\delta)\bfTh(t)^{-1};\\
  \bfA_1&=&\delta\left(\bfTh(t)^{-1}\right)^\top\left(
    \bfTh(t+\delta)^\top\bfTh(t+\delta)\frac{\partial \bff}{\partial
      \bfa}  \right.\\
  &&+ \left.
    \left(\frac{\partial \bff}{\partial \bfa}\right)^\top\bfTh(t+\delta)^\top\bfTh(t+\delta)
  \right)\bfTh(t)^{-1}.
\end{array}   
\]

Using the Taylor expansion
$\bfTh(t+\delta)=\bfTh(t)+\delta\dot{\bfTh}(t)+O(\delta^2)$ leads to
%
% XXXX
% \[
% \begin{array}{rcl}
%   \bfA_0&=&\bfI+2\delta(\dot{\bfTh}(t)\bfTh(t)^{-1})_s+O(\delta^2);\\
%   \delta \bfA_1
%   &=&\delta\left(\bfTh(t)^{-1}\right)^\top
%   \left( \bfTh(t)^\top\bfTh(t)\frac{\partial \bff}{\partial \bfa}  + 
%     \left(\frac{\partial \bff}{\partial \bfa}\right)^\top \bfTh(t)^\top\bfTh(t) 
%   \right)
%   \bfTh(t)^{-1}+O(\delta^2)\\
%   &=&\delta\left(
%     \bfTh(t)\frac{\partial \bff}{\partial \bfa}\bfTh(t)^{-1}+
%     \left(\bfTh(t)^{-1}\right)^\top\left(\frac{\partial \bff}{\partial \bfa}\right)^\top\bfTh(t)^\top
%   \right)+O(\delta^2)\\
%   &=&2\delta \left(\bfTh(t)\frac{\partial \bff}{\partial \bfa}\bfTh(t)^{-1}\right)_s+O(\delta^2).
% \end{array}   
% \]
\[
\begin{array}{rcl}
  \bfA_0&=&\bfI+2\delta(\dot{\bfTh}(t)\bfTh(t)^{-1})_s+O(\delta^2);\\
  \delta \bfA_1
  &=&\delta\left(\bfTh(t)^{-1}\right)^\top
  \left( \bfTh(t)^\top\bfTh(t)\frac{\partial \bff}{\partial \bfa}  + 
    \left(\frac{\partial \bff}{\partial \bfa}\right)^\top \bfTh(t)^\top\bfTh(t) 
  \right)
  \bfTh(t)^{-1}\\
  &&+O(\delta^2)\\
  &=&\delta\left(
    \bfTh(t)\frac{\partial \bff}{\partial \bfa}\bfTh(t)^{-1}+
    \left(\bfTh(t)^{-1}\right)^\top\left(\frac{\partial \bff}{\partial \bfa}\right)^\top\bfTh(t)^\top
  \right)+O(\delta^2)\\
  &=&2\delta \left(\bfTh(t)\frac{\partial \bff}{\partial \bfa}\bfTh(t)^{-1}\right)_s+O(\delta^2).
\end{array}   
\]
Summarizing the previous calculations, one has
\[
\bfG_k^\top\bfG_k=\bfI+2\delta \left( \left(\dot{\bfTh}(t)+
  \bfTh(t)\frac{\partial \bff}{\partial \bfa} \right)\bfTh(t)^{-1}
\right)_s+O(\delta^2),
\]
Thus, the hypothesis \textbf{(Hc1)} that $\bff$ is
contracting in the metric $\bfM$ with rate $\lambda$ implies
\[
\lambda_{\max}(\bfG_k^\top\bfG_k)\leq 1-2\delta\lambda+\epsilon(\delta),
\]
with $\lim_{\delta\to0}\frac{\epsilon(\delta)}{\delta}=0$. Letting
$\mu(\delta)=1-2\delta\lambda+\epsilon(\delta)$, one then has that
$\mu<1$ for $\delta$ sufficiently small, which in turn means that
system~(\ref{eq:SDEdiscr}) satisfies~\textbf{(Hd1)}. Applying the
discrete contraction theorem for $k=N$ leads to
%
% XXXX
% \begin{eqnarray}
%   \label{eq:res}
%   \expect_\bfx \left[d^2_{\bfM_N}(\bfa^\delta_N,\bfb^\delta_N)\right] \leq
%   \frac{2\delta C}{1-\mu(\delta)}  + 
%   \mu(\delta)^N 
%   \left[d^2_{\bfM_0}(\bfa(0),\bfb(0))-
%     \frac{2\delta C}{1-\mu(\delta)}\right]^+.
% \end{eqnarray}
\begin{eqnarray}
  \label{eq:res}
  \expect_\bfx \left[d^2_{\bfM_N}(\bfa^\delta_N,\bfb^\delta_N)\right] \leq
  \frac{2\delta C}{1-\mu(\delta)}  \nonumber\\
  +
  \mu(\delta)^N 
  \left[d^2_{\bfM_0}(\bfa(0),\bfb(0))-
    \frac{2\delta C}{1-\mu(\delta)}\right]^+.
\end{eqnarray}

On the other hand, one has, by the triangle inequality,
% XXXX
% \begin{eqnarray}
%   \expect_\bfx \left[ \dm{T}^2(\bfa(T),\bfb(T)) \right] &\leq &
%   \expect_\bfx \left[ d_{\bfM_N}^2(\bfa^\delta_N,\bfb^\delta_N)
%   \right]  \nonumber\\
%   &+&
%   \expect_\bfx \left[ d_{\bfM_N}^2(\bfa^\delta_N,\bfa(T)) \right]+
%   \expect_\bfx \left[ d_{\bfM_N}^2(\bfb^\delta_N,\bfb(T)) \right]\nonumber.
% \end{eqnarray}
\begin{eqnarray}
  \expect_\bfx \left[ \dm{T}^2(\bfa(T),\bfb(T)) \right] &\leq &
  \expect_\bfx \left[ d_{\bfM_N}^2(\bfa^\delta_N,\bfb^\delta_N)
  \right]  \nonumber\\
  &+&
  \expect_\bfx \left[ d_{\bfM_N}^2(\bfa^\delta_N,\bfa(T)) \right] \nonumber\\
  &+& 
  \expect_\bfx \left[ d_{\bfM_N}^2(\bfb^\delta_N,\bfb(T)) \right]\nonumber.
\end{eqnarray}

From equation~(\ref{eq:lim}), the second and third terms of the
right-hand side vanish when $\delta\to 0$. As for the first term,
remark that
\[
\begin{array}{lcc}
\frac{2\delta C}{1-\mu(\delta)}=\frac{2\delta
  C}{2\delta\lambda-\epsilon(\delta)}=\frac{C}{\lambda+\epsilon(\delta)/\delta}
& \xrightarrow[\delta\to 0]{} &\frac{C}{\lambda}; \ 
\\
\mu(\delta)^N=\left(1-2\delta\lambda+\epsilon(\delta)\right)^{T/\delta}
=e^{\frac{T}{\delta}\left(-2\delta\lambda+\epsilon(\delta)\right)} &
\xrightarrow[\delta\to 0]{} & e^{-2\lambda T}.
\end{array}
\]
One can thus conclude, by letting $\delta\to 0$, that
\[
\expect_\bfx \left[ \dm{T}^2(\bfa(T),\bfb(T)) \right] \leq
\frac{C}{\lambda} + e^{-2\lambda T}\left[d^2_{\bfM(0)}(\bfa(0),\bfb(0))-
    \frac{C}{\lambda}\right]^+.
\]
Integrating with respect to $\bfx$ then leads to the desired result
(\ref{eq:gen1}). Finally, (\ref{eq:gen2}) follows from (\ref{eq:gen1})
by the same calculations as in (\ref{eq:trick}) and (\ref{eq:trick2})
$\Box$

\paragraph*{Remark 3.1 [Noisy and noise-free trajectories]} If
$(\bfa,\bfb)$ represent in fact a noisy and a noise-free trajectories
then the bounds (\ref{eq:gen1}) and (\ref{eq:gen2}) are replaced by
analogous bounds where $C$ is replaced by
$C/2$~(cf. \cite{PhaX09tac}).

\paragraph*{Remark 3.2 [Relaxing the uniform bound on the noise]} As in
Remark~2.1, if the initial conditions are contained in a region $U$,
then \textbf{(Hc2)} can in fact be replaced by
% XXXX
% \[
% \forall \bfa \in U \quad \forall k\geq 0 \quad \expect\left[
%   \tr\left(\sigma(\bfa(t),t)^\top\bfM(\bfa(t),t)\sigma(\bfa(t),t)\right) \ | \
%   \bfa(0)=\bfa \right] \leq C.
% \]
\[
\forall \bfa \in U \quad \forall k\geq 0 
\]
\[ \expect\left[
  \tr\left(\sigma(\bfa(t),t)^\top\bfM(\bfa(t),t)\sigma(\bfa(t),t)\right) \ | \
  \bfa(0)=\bfa \right] \leq C.
\]

\paragraph*{Remark 3.3 [``Optimality'' of the mean square bound]}
If $\bfM$ is in fact state-independent, then the bound~(\ref{eq:gen1})
is the same as that obtained in~\cite{PhaX09tac} (cf. Theorem~2 of
that reference), which means that this bound is ``optimal'', in the
sense that it can be attained (cf. section~III-A
of~\cite{PhaX09tac}). This contrasts with the bound obtained
in~\cite{DanX13tac} (cf. Lemma~2 of that reference), which has the
same form as~(\ref{eq:gen1}) but with different constants $\lambda_1$
and $C_1$, defined -- using our notations -- as follows:
\[
\lambda_1=\lambda-\frac{\epsilon}{\beta}
\quad ; \quad
C_1=C+\frac{n\bar{m}^2\sigma^4}{2\epsilon},
\]
where $\sigma$ is a uniform upper-bound on the Frobenius norm of the
matrix $\sigma(\bfa,t)$, $\bar{m}$ is a uniform upper-bound on
$\|\bfM(\bfa,t)\|$, and $\epsilon$ is a positive constant. Note that,
for any choice of $\epsilon$, one has $\lambda_1<\lambda$ and $C_1>C$,
which yield a strictly looser bound compared to~(\ref{eq:gen1})
. Moreover, if $\epsilon$ is small, $\lambda_1$ gets closer to
$\lambda$, but $C_1$ becomes very large. On the other hand, if
$\epsilon$ is large, $C_1$ gets closer to $C$, but $\lambda_1$ becomes
very small. Thus, there is no value of $\epsilon$ for which
$\lambda_1$ and $C_1$ are arbitrarily close to $\lambda$ and $C$
respectively -- and in practice, the difference between $C_1$ and $C$
can be extremely large because of the uniform upper-bounds $\sigma$
and $\bar{m}$.

\paragraph*{Example} Following~\cite{SP12tac}, consider the following
system
\begin{equation}
  \label{eq:sys}
  \dot x_1=x_2 \sqrt{1+x_1^2} \quad ; \quad 
  \dot x_2=\frac{-x_1 x_2^2}{\sqrt{1+x_1^2}} \quad ; \quad
  y=x_1.
\end{equation}
Construct the observer
\begin{eqnarray}
  \label{eq:obshat}
  \dot{\bar{\hat{x}}}_1=\bar{\hat{x}}_2-(\bar{\hat{x}}_1-y) \quad ;
  \quad 
  \dot{\bar{\hat{x}}}_2=-(\bar{\hat{x}}_1-y),
\end{eqnarray}
\begin{eqnarray}
  \label{eq:obs}
  \hat{x}_1=\bar{\hat{x}}_1 \quad ; \quad 
  \hat{x}_2=\frac{\bar{\hat{x}}_2}{\sqrt{1+{\bar{\hat x}}_1^2}}.
\end{eqnarray}
Note that this observer differs from that of~\cite{SP12tac} : the
denominator in~(\ref{eq:obs}) is $\sqrt{1+{\bar{\hat x}}_1^2}$ instead
of $\sqrt{1+y^2}$. The observer of~\cite{SP12tac} is interesting in
that it is contracting in no state-independent metric (cf. Example 2.5
of that reference). It can be shown that this property is shared by
the modified version (\ref{eq:obshat})-(\ref{eq:obs}).

Differentiating~(\ref{eq:obs}) and replacing $\bar{\hat{x}}_1$ and
$\bar{\hat{x}}_2$ by their expressions in terms of $\hat x_1$, $\hat
x_2$, $y$, one obtains
\begin{eqnarray}
  \label{eq:obstrans}
  \dot{\hat x}_1&=&\dot{\bar{\hat x}}_1
  =\bar{\hat x}_2-(\bar{\hat{x}}_1-y)
  =\hat x_2 \sqrt{1+\hat{x}_1^2} - (\hat{x}_1-y)\quad ;\nonumber\\
  \dot{\hat x}_2&=&\frac{\dot{\bar{\hat{x}}}_2}{\sqrt{1+\hat{x}_1^2}} 
  - \frac{\bar{\hat x}_2 \hat{x}_1\dot{\hat
      x}_1}{(1+\hat{x}_1^2)^{3/2}}\\
  &=&
  -\frac{(\hat{x}_1-y)\left(\hat{x}_1\hat{x}_2-\sqrt{1+\hat{x}_1^2}\right)}{1+\hat{x}_1^2} 
  - \frac{\hat{x}_1\hat{x}_2^2}{\sqrt{1+\hat{x}_1^2}}\nonumber.
\end{eqnarray}

Observe that $(x_1,x_2)$ is a particular solution
of~(\ref{eq:obstrans}). To show the contraction behavior
of~(\ref{eq:obstrans}), consider the following nonlinear transform
\begin{eqnarray}
  \label{eq:check}
  \check{x}_1&=&-3\hat{x}_1+5\hat{x}_2\sqrt{1+\hat{x}_1^2},\nonumber\\
  \check{x}_2&=&3\hat{x}_1+2\hat{x}_2\sqrt{1+\hat{x}_1^2}.
\end{eqnarray}
From~(\ref{eq:obs}), one has 
\[
(\check{x}_1,\check{x}_2)^\top=\bfP \cdot (\bar{\hat x}_1,\bar{\hat x}_2)^\top,
\]
where $\bfP$ is the $2\times 2$ constant matrix 
$\left(\begin{array}{cc}
    -3&5\\
    3&2
\end{array}\right)$. Thus
\[
(\dot{\check x}_1,\dot{\check x}_2)^\top
=\bfP\cdot (\dot{\bar{\hat x}}_1,\dot{\bar{\hat x}}_2)^\top
=\bfP\bfQ \cdot (\bar{\hat x}_1,\bar{\hat x}_2)^\top
=\bfP\bfQ\bfP^{-1} \cdot (\check{x}_1,\check{x}_2)^\top,
\]
where the second inequality comes from~(\ref{eq:obstrans}) with
$\bfQ=\left(\begin{array}{cc}
    -1&1\\
    -1&0
\end{array}\right)$.
A numerical computation shows that the eigenvalues of the symmetric
part of $\bfP\bfQ\bfP^{-1}$ are $(-0.24,-0.76)$, which means that
system~$(\check{x}_1,\check{x}_2)$ is contracting with rate $0.24$ in
the identity metric. From~(\ref{eq:obs}), one finally has that
system~$(\hat{x}_1,\hat{x}_2)$ is contracting with rate $0.24$ in the
metric
\[
\bfM=\bfTh^\top\bfP^\top\bfP\bfTh,\quad \textrm{where}\ 
\bfTh=\left(\begin{array}{cc}
    1&0\\
    \frac{-\hat{x}_1\hat{x}_2}{\sqrt{1+\hat{x}_1^2}}&\sqrt{1+\hat{x}_1^2}
  \end{array}\right).
\]

Let us now study the convergence properties of the observer when the
measure $y_p$ is \emph{corrupted by white noise} as $y_p=y+S\xi$,
where $y=x_1$ is the unperturbed measure, $\xi$ is a ``white noise''
of variance 1 and $S$ is the noise intensity. Using the formal rule
$\ud W=\xi\ud t$, equations~(\ref{eq:obshat}) are transformed into
% XXXX
% \begin{eqnarray}
%   \label{eq:obshatnoise}
%   \ud \bar{\hat{x}}_1=(\bar{\hat{x}}_2-(\bar{\hat{x}}_1-y))\ud t+S\ud W \quad ;
%   \quad 
%   \ud \bar{\hat{x}}_2=-(\bar{\hat{x}}_1-y)\ud t+S\ud W.
% \end{eqnarray}
\begin{eqnarray}
  \label{eq:obshatnoise}
  \ud \bar{\hat{x}}_1&=&(\bar{\hat{x}}_2-(\bar{\hat{x}}_1-y))\ud t+S\ud W \nonumber\\
  \ud \bar{\hat{x}}_2&=&-(\bar{\hat{x}}_1-y)\ud t+S\ud W.
\end{eqnarray}

The observer equations~(\ref{eq:obstrans}) become
%XXXX
% \begin{eqnarray}
%   \label{eq:obstransnoise}
%   \ud \hat x_1&=&\left[\hat x_2 \sqrt{1+\hat{x}_1^2} -
%     (\hat{x}_1-y)\right]\ud t + S\ud W\quad ;\nonumber\\
%   \ud \hat x_2&=&-\left[\frac{(\hat{x}_1-y)\left(\hat{x}_1\hat{x}_2-\sqrt{1+\hat{x}_1^2}\right)}{1+\hat{x}_1^2} 
%     + \frac{\hat{x}_1\hat{x}_2^2}{\sqrt{1+\hat{x}_1^2}}\right]\ud t
%   -S\left[\frac{\hat{x}_1\hat{x}_2-\sqrt{1+\hat{x}_1^2}}{1+\hat{x}_1^2}\right]\ud W
% \nonumber.
% \end{eqnarray}
\begin{eqnarray}
  \label{eq:obstransnoise}
  \ud \hat x_1&=&\left[\hat x_2 \sqrt{1+\hat{x}_1^2} -
    (\hat{x}_1-y)\right]\ud t + S\ud W\quad ;\nonumber\\
  \ud \hat x_2&=&-\left[\frac{(\hat{x}_1-y)\left(\hat{x}_1\hat{x}_2-\sqrt{1+\hat{x}_1^2}\right)}{1+\hat{x}_1^2} 
    + \frac{\hat{x}_1\hat{x}_2^2}{\sqrt{1+\hat{x}_1^2}}\right]\ud t\nonumber\\
  &-&S\left[\frac{\hat{x}_1\hat{x}_2-\sqrt{1+\hat{x}_1^2}}{1+\hat{x}_1^2}\right]\ud W
\nonumber.
\end{eqnarray}

One is now in the settings of Theorem~\ref{theo:main-cont} with
\[
\sigma(\hat x_1,\hat x_2)=
\left(S,S\frac{\hat{x}_1\hat{x}_2-\sqrt{1+\hat{x}_1^2}}{1+\hat{x}_1^2}\right)^\top.
\]
From the above expression, it can be shown algebraically that
$\sup_{a,b}\sigma(a,b)^\top\bfM(a,b)\sigma(a,b)=15.2S^2$. 

\begin{figure}[htp]
    \centering
    \hspace{0cm}\textbf{A}\hspace{3cm}\textbf{B}\\
    \vspace{0.1cm}
    \includegraphics[width=3.5cm]{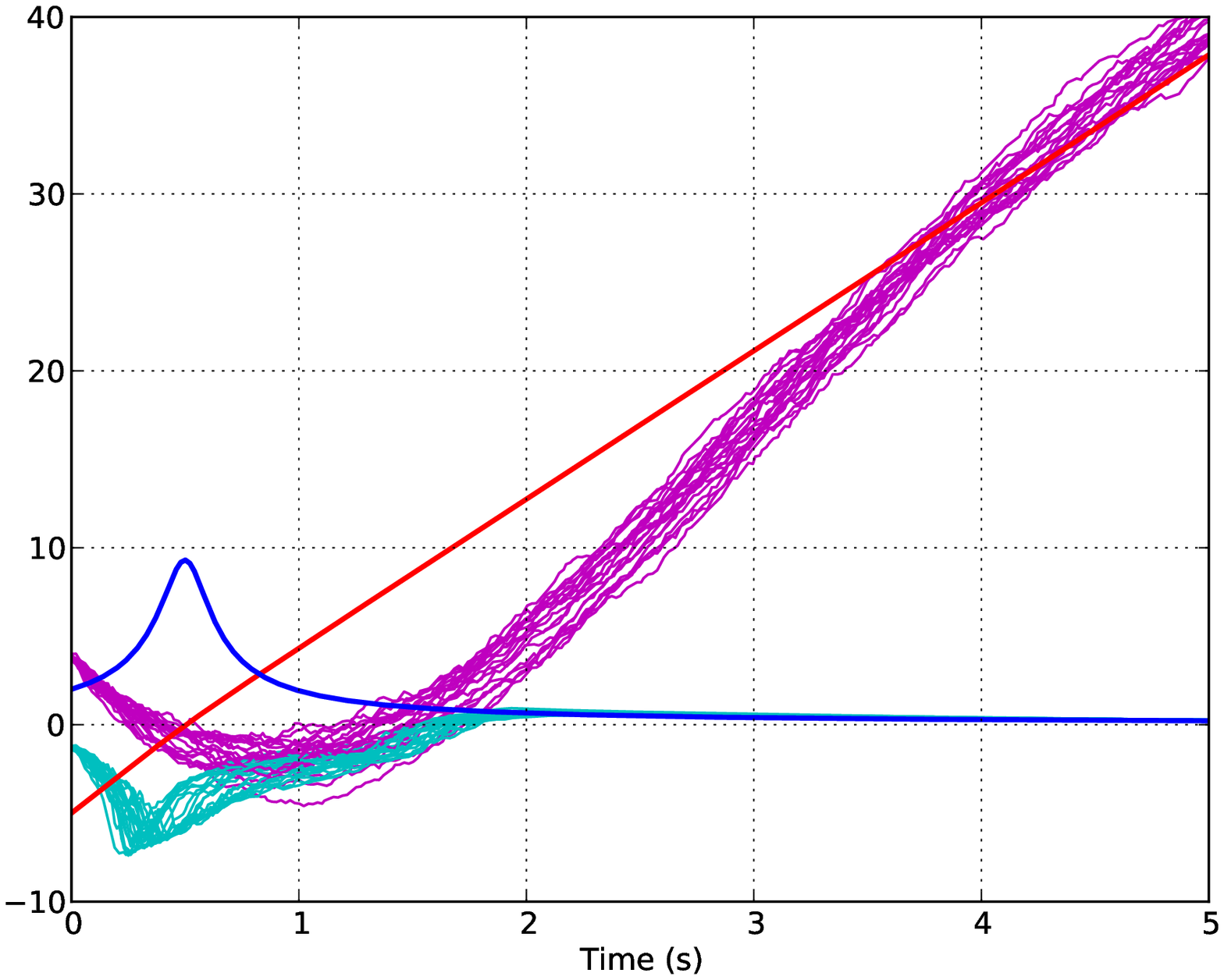}
    %\hspace{2cm}
    \includegraphics[width=3.5cm]{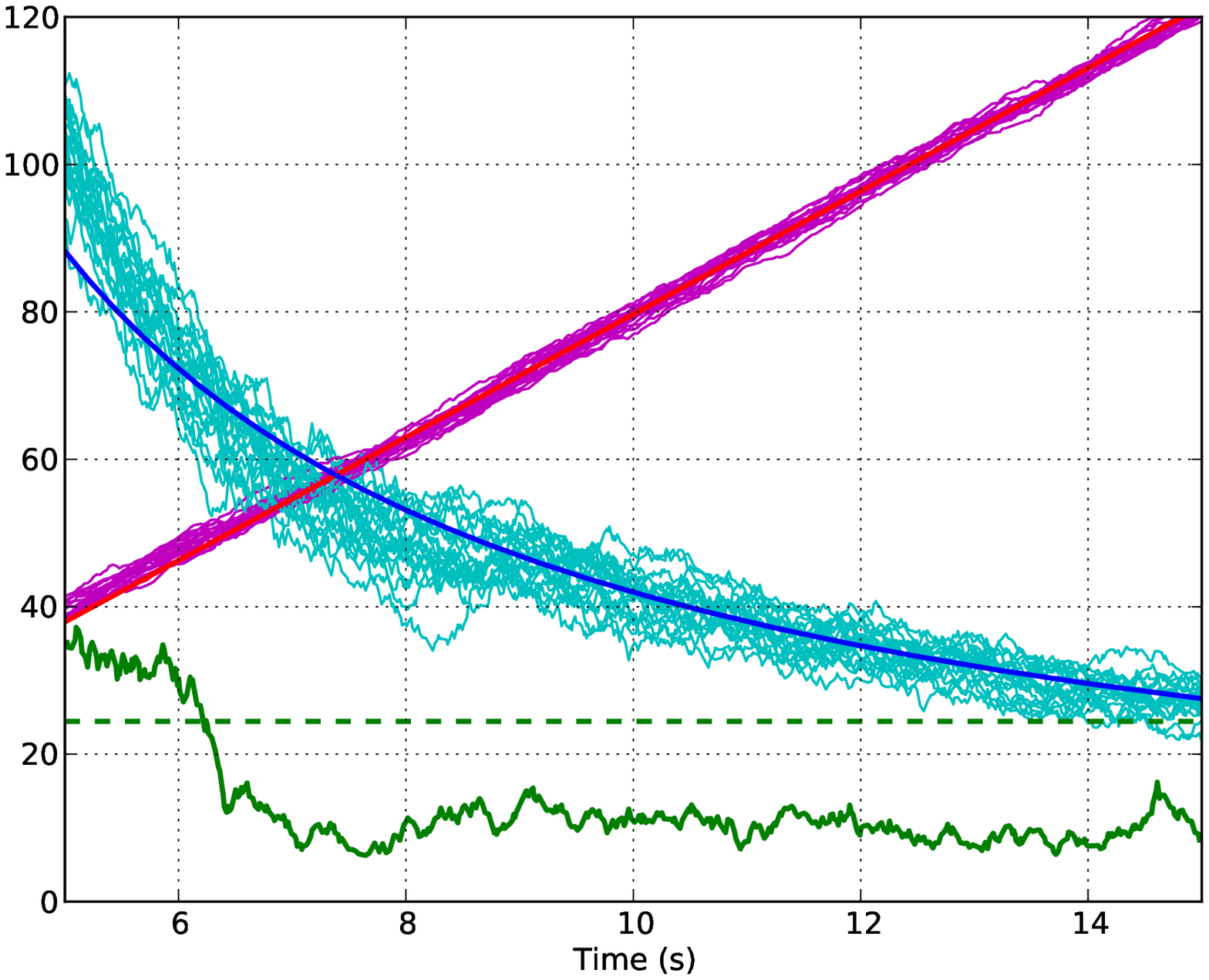}\\
    \caption{Simulations for the observer studied in the
      text. \textbf{A}: evolution of the systems for
      $t\in[0$\,s$,5$\,s$]$.  Equations~(\ref{eq:sys}) were integrated
      using the Euler method with time step $\Delta t=0.01$\,s (red
      line: $x_1$; blue line: $x_2$). Equations~(\ref{eq:obshatnoise})
      were integrated using the Euler-Maruyama
      scheme~(cf. \cite{KP92book}) with the same time step $\Delta
      t=0.01$\,s. We plotted 20 sample trajectories for noise
      intensity $S=1$ starting from the same deterministic initial
      values $(\hat{x}_1(0),\hat{x}_2(0))$ (magenta lines:
      $\hat{x}_1$; cyan lines: $\hat{x}_2$).  \textbf{B}: evolution of
      the systems for $t\in[5$\,s$,15$\,s$]$. Note that, for clarity,
      the values of $x_2$ and $\hat{x}_2$ were multiplied by 400 in
      this plot. To assess the theoretical bounds, we plotted the
      sample mean square error $(x_1-\hat{x}_1)^2+(x_2-\hat{x}_2)^2$
      (plain green line) and the theoretical bound after transients
      given by equation~(\ref{eq:ult}) (dashed green line). For
      clarity, these values were multiplied by 10.}
    \label{fig:simu}
\end{figure}

We now make the assumption that $\|\hat x_2\|$ is uniformly upper-bounded
by a constant $B$ (which can indeed be shown using an independent
method, see also simulations in Fig.~\ref{fig:simu}). Then, it can be
shown that, uniformly,
\[
\|\bfTh^\top\bfP^\top\bfP\bfTh\bfx\|^2 \geq \gamma(B)\|\bfx\|^2.
\]

One thus can apply Theorem~\ref{theo:main-cont} and obtain the
bound~(\ref{eq:gen2}) with $\lambda=0.24$, $C=15.2S^2$ and
$\beta=\gamma(B)$. Note that, for $t\to\infty$, one has $\hat x_2\to
0$, such that one has the bound $B=0$, which in turn corresponds to
$\gamma(B)=12.95$. The bound after exponential transients is then
given by (cf. Fig.~\ref{fig:simu} for numerical simulations)
\begin{equation}
  \label{eq:ult}
  \frac{C}{2\beta\lambda}=2.45S^2.
\end{equation}

\section{Conclusion}
\label{sec:conclusion}

We have established the stochastic contraction theorems in the case of
general time- and state-dependent Riemannian metrics. In the limit
when the metric becomes linear (state-independent), the bounds we
derived are the same as those obtained in~\cite{PhaX09tac}, which
means that they are ``\emph{optimal}'', in the sense that they can be
attained.  This development allows extending the applicability of
contraction analysis to a significantly wider range of nonlinear
stochastic dynamics, such as stochastic observers or networks of noisy
nonlinear oscillators.

\section*{Acknowledgments} We thank J.-M. Mirebeau for his help with
the proof of Proposition~2. QCP was supported by a JSPS postdoctoral
fellowship.

{\small
\bibliographystyle{abbrv}
\bibliography{../../ynl}
}

\end{document}